\documentclass[12pt]{amsart}
\usepackage{amsfonts}
\usepackage{amssymb}
\usepackage{amsmath}
\usepackage[all]{xy}
\usepackage[utf8]{inputenc}
\usepackage[T1]{fontenc}
\usepackage[english]{babel}
\usepackage{xspace}
\usepackage{times}
\usepackage{amsthm}
\usepackage{amscd}
\usepackage{amsbsy}
\usepackage{mathrsfs}
\usepackage{relsize}
\usepackage[center,small,sc]{caption}
\usepackage{pgf}
\usepackage{lmodern}
\usepackage{xcolor}

\usepackage{fancyhdr}
\usepackage{abstract}

\newtheorem{th1}{{\bf Theorem}}[section]

\newtheorem{le1}[th1]{{\bf Lemma}}
\newtheorem{pr1}[th1]{{\bf Proposition}}

\theoremstyle{remark}

\theoremstyle{definition}
\newtheorem{de1}[th1]{Definition}

\newtheorem{th2}{{\bf Theorem}}[subsection]

\newtheorem{pr2}[th2]{{\bf Proposition}}

\theoremstyle{remark}

\theoremstyle{definition}

\newtheorem{th3}{{\bf Theorem}}[subsubsection]
\newtheorem{le3}[th3]{{\bf Lemma}}
\newtheorem{pr3}[th3]{{\bf Proposition}}

\theoremstyle{remark}
\newtheorem{re3}[th3]{Remark}

\theoremstyle{definition}

\newtheorem{No3}[th3]{Notations}

\DeclareGraphicsExtensions{.jpg,.mps,.pdf,.png,.gif}

\begin{document}

\title[Equivariant cohomology and depth]{Equivariant cohomology and depth}
\author{Dorra BOURGUIBA and Said ZARATI}
\address{Université Tunis-ElManar, Facult\'e des Sciences de Tunis,
D\'epartement de Math\'ematiques.
TN-2092 Tunis, TUNISIE.}
\email{dorra.bourguiba@fst.utm.tn}
\email{said.zarati@fst.rnu.tn}
\dedicatory{{\it In memory of our colleague Sami Hammouda}}

\subjclass{Algebraic topology}
\date{}
\maketitle
\begin{abstract} { \it Let $n \geq 1$ be an integer, let $V=(\mathbb{Z}/2\mathbb{Z})^{n}$ and let $X$ be a $V$-CW-complex. If $X$ is a finite $CW$-complexe, the equivariant modulo $2$ cohomology of the  $V$-CW-complexe $X$, denoted by $H_{V}^{*}(X, \mathbb{F}_{2})$, is a finite type module over the modulo $2$ cohomology of the group $V$, denoted by $H^{*}(V, \mathbb{F}_{2})$. Let $dth_{H^{*}V}H_{V}^{*}(X, \mathbb{F}_{2})$ be the depth of the finite type $H^{*}(V, \mathbb{F}_{2})$-module $H_{V}^{*}(X, \mathbb{F}_{2})$ relatively to the augmentation ideal, $\widetilde{H^{*}}(V, \mathbb{F}_{2})$, of $H^*(V, \mathbb{F}_{2})$.
\medskip\\
Let $W$ be a subgroup of $V$, the finite $V$-CW-complex $X$ is a finite $W$-CW-complex, we denote by \\ $dth_{H^{*}W}H_{W}^{*}(X, \mathbb{F}_{2})$ the depth of the finite type $H^{*}(W, \mathbb{F}_{2})$-module  $H_{W}^{*}(X, \mathbb{F}_{2})$ relatively to the augmentation ideal , $\widetilde{H^{*}}(W, \mathbb{F}_{2})$, of $H^*(W, \mathbb{F}_{2})$.
\medskip\\
The aim of this paper is to study a question raised in $2005$ in \cite{BouHZ} concerning the comparaison of  $dth_{H^{*}V} H_{V}^{*}(X, \mathbb{F}_{2})$ and $dth_{H^{*}W}H_{W}^{*}(X, \mathbb{F}_{2})$.
\medskip\\
We prove the following result:
\medskip\\
{\bf Theorem}: For every subgroup $W$ of $V$, we have: $dth_{H^{*}W}H_{W}^{*}(X, \mathbb{F}_{2}) \leq dth_{H^{*}V} H_{V}^{*}(X, \mathbb{F}_{2})  $.}
\end{abstract}
\section{Introduction.}
Let $V$ be an elementary abelian 2-{group}, that is $V \cong (Z/2Z)^n$, for an integer $n \geq 1$. This integer $n$ is called the rank of the group $V$, and denoted $rk(V)$, or the dimension of the $\mathbb{F}_{2}$-vector space $V$ and denoted by $\dim_{\mathbb{F}_{2}}V$. Let $EV$ be a contactible space on which the group $V$ acts freely and let $BV = EV/V$ be a classifying space for the group $V$. Recall that the modulo~2 cohomology of the group $V$, $H^*V = H^*(BV; \mathbb{F}_{2})$, is isomorphic to a graded polynomial algebra $\mathbb{F}_{2}[t_1, \dots, t_n]$ on generators $t_i$, $1 \le i \le n$, of degree 1.\\
Let $X$ be a $V$-CW-complex (for a reference of the notion of equivariant CW-complex, see for example [\cite{TD}, Chap. II, Sect. 1]); we denote by $X_{hV}  = EV \times_{V}X$ the quotient of $EV \times X$ by the diagonal action of the group $V$. The modulo~2 cohomology of the space $X_{hV}$, denoted by $H^{*}_{V}X=H_{V}^{*}(X, \mathbb{F}_{2})$ and called modulo~2 equivariant cohomology of $X$, is naturally an $\mathrm{H}^{*}V$-module via the morphism induced, in modulo~2 cohomology, by the projection $\pi_{V}: X_{hV} \rightarrow \{*\}_{hV}=BV$.
\medskip\\
If $X$ is a finite $V$-CW-complex, we verify that $H^{*}_{V}X$ is an $\mathrm{H}^{*}V$-module of finite type (see \cite{Q}); we then denote by $dth_{H^*V}H^*_{V}X$ the depth of $H^{*}_{V}X$ relatively to the augmentation ideal, $\widetilde{H^{*}}V$, of $\mathrm{H}^{*}V$ (see \cite{Ma}).\\
Let $W$ be a subgroup of $V$, the finite $V$-CW-complex $X$ is also a finite $W$-CW-complex; so we have two integers namely: $dth_{H^*V}H^*_{V}X$ and $dth_{H^*W}H^*_{W}X$. The comparaison of these two integers is a natural question asked in $2005$ in \cite{BouHZ} where the authors conjectured the following inequality:
\begin{center}
Let $V$ be an elementary abelian $2$-group and $X$ be a finite $V$-CW-complex. Then, for every subgroup $W$ of $V$, we have:
 $dth_{H^*W} H^{*}_{W}X \leq dth_{H^*V} H^{*}_{V}X $.
\end{center}
\smallskip
The purpose of this paper is to prove the previous inequality and to deduce some consequences. Our method is mainly to replace $\mathrm{H}^{*}V$ with the Dickson algebra $DV$ (see 2.3.2). For this we use the results of Serre \cite{Se1}, \cite{Se 2}  and Bourguiba-Zarati \cite{BouZ}.
\smallskip\\
The paper is structured as follows. In the first paragraph we recall the notion of depth of a finite type $H^*V$-module, give its main properties and fix some notations. The second paragraph is  devoted to some technical results. In the third paragraph we give the proof of the main result. The fourth paragraph is dedicated to some remarks and consequences of the main result.
\bigskip\\
\textbf{Acknowledgements} The first  author is a member of  the "Laboratoire de Recherche {\bf LATAO} code LR11ES12" of the Faculty of Sciences of Tunis at the  University of Tunis El-Manar.\\
The authors would like to thank Professeur Jean Lannes for discussions they had on this subject during the annual colloquium "CSMT2022" of the Tunisian Mathematical Society.
\section{Definitions, notations and recalls }
Let $V$ be an elementary abelian $2$-group. We recall in this section the definitions and some properties of regular sequences and depth of an $\mathrm{H}^{*}V$-module  (for clarity we only limit to $\mathrm{H}^{*}V$-modules).
\subsection{Depth of a finite type $\mathrm{H}^{*}V$-module}
Let $M$ be a finite type $\mathrm{H}^{*}V$-module. A sequence of homogeneous elements $\{\alpha_1,\ldots,\alpha_k \}$ of stricly positive degrees of $\mathrm{H}^{*}V$ is said to be regular for $M$ if:
\begin{itemize}
\item[--] $\alpha_1$ is not a zero divisor on $M$,
\item[--] for each $i = 2,\ldots,k$, $\alpha_i$ is not a zero divisor on $M/(\alpha_1,\ldots,\alpha_{i-1})M$, where ($\alpha_1,\ldots, \alpha_{i-1})$ is the ideal of $\mathrm{H}^{*}V$ generated by $\alpha_1,\ldots,\alpha_{i-1}$.
\end{itemize}
\smallskip
The integer $k$ is called the length of the sequence $\{\alpha_1,\ldots, \alpha_k \}$.
\smallskip\\
A regular sequence $\{\alpha_1,\ldots,\alpha_k \}$ of $M$ is said to be maximal if, for any $\alpha \in \widetilde{H^{*}}V$, the sequence $\{\alpha_1,\ldots,\alpha_k, \alpha\}$ is no longer regular for $M$.
\smallskip\\
Let $M$ be an $\mathrm{H}^{*}V$-module of finite type; since $\mathrm{H}^{*}V$ is a graded commutative noetherian algebra so any maximal regular sequence for $M$ have the same finite length, (see \cite{Ma}).  This common length is called depth of $M$ and denoted by $dth_{H^{*}V} M$. By convention,  the depth of the trivial $\mathrm{H}^{*}V$-module is $+ \infty$.
\subsection{Some properties of depth}
Recall that any $\mathrm{H}^{*}V$-module $M$ has a projective dimension, denoted by $dp_{H^{*}V}M$, and defined as the minimal length of a projective $\mathrm{H}^{*}V$-resolution of $M$ (see for example \cite{Bo}). The depth of an $\mathrm{H}^{*}V$-module of finite type is related to the projective dimension by the Auslander-Buschbaum formula (see \cite{Ma}):
\smallskip
\begin{center}
 $(2.2.1) \quad  dth_{H^{*}V}M + dp_{H^{*}V}M = rk(V)$
\end{center}
\smallskip
By convention, the projective dimension of the trivial $\mathrm{H}^{*}V$-module is $-\infty$. An $\mathrm{H}^{*}V$-module of finite type $M$ is said Cohen-Macaulay if $dp_{H^{*}V} M = 0$ (or equivalently, if $dth_{H^{*}V}M = rk(V)$).
\smallskip\\
Let's denote $Ext_{{H^{*}V}}^{i}(-, M)$ \big(resp. $Tor_{_{i}}^{H^{*}V}(-, M))$, $i \ge 0$ \big), the derived functors of the functor
 $Hom_{H{^*}V}(-, M)$  \big(resp. of functor $- \otimes_{H{^*}V}M$ \big). One characterization of the depth of a finite type $\mathrm{H}^{*}V$-module $M$ is given in \cite{Ma}, in terms of $Ext$, by:
\smallskip
\begin{center}
 $(2.2.2) \quad dth_{H^{*}V}M = \inf \{i\; \mid \; Ext_{H^{*}V}^{i}(\mathbb{F}_{2}, M) \neq 0\}.$
\end{center}
\smallskip
Let $0 \rightarrow M' \rightarrow M \rightarrow M" \rightarrow 0$ be an exact sequence of $\mathrm{H}^{*}V$-modules of finite type. As a consequence of the previous characterization, we have:
\smallskip
\begin{center}
 $(2.2.3) \quad   dth_{H^{*}V}M \geq \inf \big(dth_{H^{*}V} M',\; dth_{H^{*}V} M" \big)$
\end{center}
\smallskip
Another characterization of the depth of an $\mathrm{H}^{*}V$-module (or equivalently of the  projective dimension of an $\mathrm{H}^{*}V$-module) is given in \cite{Bo}, in terms of $Tor$, by:
\smallskip
\begin{center}
$(2.2.4)\quad (dp_{H^{*}V}M \le k) \Leftrightarrow \big(Tor_{i}^{H^{*}V}(\mathbb{F}_{2}, M) = 0 \; \mbox{pour} \; i \ge k+1 \big)$
\end{center}
\subsection{On the computation of the depth of an $\mathrm{H}^{*}V$-module of finite type}
Computing the depth of an $\mathrm{H}^{*}V$-module of finite type is a rather difficult problem (see \cite{B}, \cite{Be}, \cite{BouZ}, \cite{Du}, \cite{BrH}, \cite{Sm}); however, we have some results.
\subsubsection{\bf Depth and tensor product}
Let $V$ be an elementary abelian 2-group, $W$ a subgroup of $V$ and $M$ an $\mathrm{H}^{*}V$-module. Note that any $\mathrm{H}^{*}V$-module $M$ is an $\mathrm{H}^{*}(V/W)$-module  via the morphism $\pi^* :  H{^*}(V/W) \rightarrow H^*V$ induced, on cohomology, by the canonical projection $\pi: V\rightarrow V/W$. We denote:
\[
   \mathcal{T}^{H{^*}(V/W)}(M ) = \{x \in M \; \mid \; \alpha x = 0,\;  \forall \alpha \in \widetilde{H^{*}}(V/W)\}
\]
We check that $\mathcal{T}^{H{^*}(V/W)}(M )$ is a sub-$H^*V$-module of $M$, called trivial part of $M$ ($M$ considered as an
$\mathrm{H}^{*}(V/W)$-module). The module $\mathcal{T}^{H{^*}(V/W)}(M )$ is also an $\mathrm{H}^{*}W$-module.\\
It is clear that if $M$ is an $\mathrm{H}^{*}V$-module of finite type,  then $\mathcal{T}^{H{^*}(V/W)}(M )$ is an $\mathrm{H}^{*}W$-module of finite type.
\smallskip\\
We have the following result which compare, in some cases, the depth of an $\mathrm{H}^{*}V$-module of finite type $M$ and the depth of the
$\mathrm{H}^{*}W$-module  of finite type $H{^*}W \otimes_{H{^*}V}M \cong \mathbb{F}_{2} \otimes_{H{^*}(V/W)} M$.
\begin{pr3} (\cite{BouHZ}) Let $V$ be an elementary abelian $2$-group, let $W$ be a subgroup of $V$ of codimension one and let $M$ be an $\mathrm{H}^{*}V$-module of finite type such that $\mathcal{T}^{H{^*}(V/W)}(M ) = 0$. Then, we have: $dth_{H^{*}V} M = dth_{H^{*}W} \big(H{^*}W \otimes_{H{^*}V}M \big) + 1$.
\medskip\\
In particular, if $M$ is Cohen-Macaulay, then so is $H{^*}W \otimes_{H{^*}V}M$.
\end{pr3}
\begin{re3}
The condition $\mathcal{T}^{H{^*}(V/W)}(M ) = 0$ in proposition 2.3.1.1 is necessary since, for example, the
proposition 2.3.1.1 is no longer true for the finite $\mathrm{H}^{*}V$-modules.
\end{re3}
\begin{le3} (\cite{BouZ})  Let $V$ be an elementary abelian $2$-group, let $W$ be a subgroup of $V$ and let $N$ be an $\mathrm{H}^{*}(V/W)$-module of finite type. Then, $$dth_{H^{*}V} (H{^*}V \otimes_{H{^*}(V/W)}N) =  dth_{H^{*}(V/W)} N + rg(W)$$
\end{le3}
\subsubsection{\bf Depth and Dickson algebra}
Let $V$ be an elementary abelian 2-group of rank $n \geq 1$ and $GL(V)$ the automorphism group of $V$. The group $GL(V)$ acts on $\mathrm{H}^{*}V$ and the ring of invariants $DV=(H^*V)^{GL(V)}$, called Dickson algebra, is a graded polynomial ring $DV=\mathbb{F}_{2}[c_1,...,c_n]$, on Dickson classes $c_{i},\; 1 \leq i \leq n=rg(V)$, of degree $2^{n} - 2^{n - i}$. In particular, $c_{n} = \prod_{u \in H^{1}V \ \{0\}}u$. In certain cases, we will specify, not to be confused, the notation $c_{i}$ by $c_{V, i}$ (see \cite{Di}, \cite{Mu}).
\smallskip\\
It's known that $\mathrm{H}^{*}V$ is $DV$-module via the natural inclusion $DV \hookrightarrow H^*V$ and by, \cite{Se1}, \cite{R} (see also  \cite{St}, \cite{HLS}, \cite{BrZ}),  $\mathrm{H}^{*}V$ is a free $DV$-module of finite type. It follows that, if $M$ is an $\mathrm{H}^{*}V$-module of finite type, then $M$ is a $DV$-module of finite type and we have, by \cite{Se 2}, $dth_{H^{*}V}M = dth_{DV}M$.
\medskip\\
Let $A$ be the modulo $2$ Steenrod algebra and $\mathcal{U}$ the category of unstable $A$-modules (see \cite{LZ2}, \cite{Sc}). We denote $V$-$\mathcal{U}$ the category of unstable $\mathrm{H}^{*}V$-$\mathrm{A}$-modules and $V_{tf}$-$\mathcal{U}$ the full sub-category of $V$-$\mathcal{U}$ whose objects are the unstable $\mathrm{H}^{*}V$-$\mathrm{A}$-modules which are of finite type as $H^*V$-modules (see \cite{LZ1}, \cite{HLS}).
For example, if $X$ is a finite $V$-CW-complex then, $H^{*}_{V}X$ is an object of the category $V_{tf}$-$\mathcal{U}$.
In \cite{BouZ} (see theorem 1.3 of \cite{B}), we prove :
\medskip
\begin{th3} Let $V$ be an elementary abelian $2$-group and $M$ is an unstable $\mathrm{H}^{*}V$-$\mathrm{A}$-module which is of finite type as $\mathrm{H}^{*}V$-module. Then, the following statements are equivalent:
\begin{itemize}
\item[(i)] $dth_{H^{*}V} M \geq k$.
\medskip
\item[(ii)]the sequence $\{c_1,\ldots,c_k \}$ is regular for $M$.
\end{itemize}
\end{th3}
\medskip
Let $M$ be an unstable $\mathrm{H}^{*}V$-$\mathrm{A}$-module which is of finite type as $\mathrm{H}^{*}V$-module and such that \\$dth_{H^{*}V} M \geq k$. the following result describes $M$, up to isomorphism, as $DV$-module.
\medskip\\
 Before stating the result, we need to introduce some notations.
\smallskip
\begin{No3}
Let $\mathbb{F}_{2}[\xi_{1},..., \xi_{s}]$ be a graded polynomial algebra over the field $\mathbb{F}_{2}$ and $E$ an $\mathbb{F}_{2}[\xi_{1},..., \xi_{s}]$-module. To simplify the notations, we denote:
$$\overline{E}^{(\xi_{1},..., \xi_{i})} = E / \widetilde{\mathbb{F}}_{2}[\xi_{1},..., \xi_{i}].E, \; 1 \leq i \leq s$$
where $\widetilde{\mathbb{F}}_{2}[\xi_{1},..., \xi_{i}]$ is the augmentation ideal of $\mathbb{F}_{2}[\xi_{1},..., \xi_{i}]$.
\medskip\\
We verify that:
$$\begin{array}{lll}
                  \overline{E}^{(\xi_{1},..., \xi_{i})} & = & \mathbb{F}_{2} {\underset{\mathbb{F}_{2}[\xi_{1},..., \xi_{i}]}{\otimes}}E = \mathbb{F}_{2}[\xi_{i+1},..., \xi_{s}] {\underset{\mathbb{F}_{2}[\xi_{1},..., \xi_{s}]}{\otimes}}E
\medskip\\
&=& Tor_{0}^{\mathbb{F}_{2}[\xi_{1},..., \xi_{s}]}(\mathbb{F}_{2}[\xi_{i+1},..., \xi_{s}],\; E) = Tor_{0}^{\mathbb{F}_{2}[\xi_{1},..., \xi_{i}]}(\mathbb{F}_{2},\; E)
 \end{array}$$
 We also denote :
$$\begin{array}{lll}
                   \mathcal{T}^{\xi_{i}}(E) &=& Ker( \xi_{i}.: E \rightarrow E, e \mapsto \xi_{i}.e )
\medskip\\
 \mathcal{T}^{(\xi_{1},..., \xi_{i})}(E)& = & \mathcal{T}^{\xi_{1}}\big( \mathcal{T}^{(\xi_{2},..., \xi_{i})}(E)\big)
 \end{array}$$
We verify that:
$$\Sigma^{\mid \xi_{i} \mid } \mathcal{T}^{\xi_{i}}(E) = Tor_{1}^{\mathbb{F}_{2}[\xi_{i}]}(\mathbb{F}_{2}, E)$$
$\Sigma$ is the suspension functor and $\mid \xi_{i} \mid$ the degree of $ \xi_{i}$.
\end{No3}
\smallskip
\begin{le3} Let $V$ be an elementary abelian $2$-group of rank $n \geq1$, let $k \in \mathbb{N}, 1 \leq k \leq n$, and let $M$ be an unstable $\mathrm{H}^{*}V$-$\mathrm{A}$-module which is of finite type as $\mathrm{H}^{*}V$-module. Then the following properties are equivalent:
\begin{itemize}
\item[(i)] $dth_{H^{*}V}M \geq k$.
\medskip
\item[(ii)] The module $M$ is isomorphic, as a $DV$-module, to: $\mathbb{F}_{2}[c_1,...,c_k] \otimes \overline{M}^{(c_1,...,c_k)}$.
\end{itemize}
\end{le3}
\begin{proof} The implication $(ii) \Rightarrow (i)$ is clear. The implication $(i) \Rightarrow (ii)$ will be proved by induction on $k$.
\medskip\\
{\bf The case $k=1$}. Suppose that $dth_{H^{*}V}M = 1$, since $M$ is an object of the category $V_{tf}$-$\mathcal{U}$, by \cite{BouZ}, the sequence $\{c_{1}\}$ is regular for $M$; it means that the multiplication by $c_{1}$ is injective in $M$. Consider the beginning of a projective $DV$-resolution of $M$ (recall that the Dickson algebra $DV \cong \mathbb{F}_{2}[c_1,...,c_n]$).
\begin{center}
$(\mathcal{R}): \; \xymatrix{ Ker(\varphi) \, \ar@{^{(}->}[r]& \mathbb{F}_{2}[c_1,...,c_n] \otimes N \ar@{->>}[rr]^-{\varphi}&& M  }$
\end{center}
where $N$ is a graded $\mathbb{F}_{2}$-vector space of finite dimension.
\medskip\\
By considering "$\mathbb{F}_{2} \otimes_{\mathbb{F}_{2}[c_{1}]}(\mathcal{R})$", and since the multiplication by $c_{1}$ is injective in $M$, we obtain this exact sequence:
 \begin{center}
$(\overline{\mathcal{R}}): \; \xymatrix{ \overline{Ker(\varphi)}^{(c_{1})} \, \ar@{^{(}->}[r]& \mathbb{F}_{2}[c_2,...,c_n] \otimes N \ar@{->>}[rr]^-{\overline{\varphi}}&& \overline{M }^{(c_{1})} }$
\end{center}
and the following commutative diagram whose bottom line is "$\mathbb{F}_{2}[c_{1}] \otimes (\overline{\mathcal{R}})$".
\begin{center}
$\xymatrix{ Ker(\varphi) \, \ar@{^{(}->}[r]  \, \ar@{^{(}->}[d]^{ \alpha}&  \mathbb{F}_{2}[c_1,...,c_n] \otimes N \ar[d]^{=} \ar@{->>}[rr]^-{\varphi}&& M \ar@{->>}[d]^{\beta} \\
\mathbb{F}_{2}[c_{1}] \otimes \overline{Ker(\varphi)}^{(c_{1})} \, \ar@{^{(}->}[r]& F _2[c_1,...,c_n] \otimes N \ar@{->>}[rr]&& \mathbb{F}_{2}[c_{1}] \otimes \overline{M }^{(c_{1})} }$
\end{center}
\smallskip
Again because of the injectivity of the multiplication by $c_{1}$ in $M$, we verify that the homomorphism of $DV$-modules $\beta: M \rightarrow \mathbb{F}_{2}[c_{1}] \otimes \overline{M }^{(c_{1})}$ is an isomorphism.
\medskip\\
Assume that the implication $(i) \Rightarrow (ii)$ is true for $k$ and let us show it for $k+1$. Let $M$ be an unstable $\mathrm{H}^{*}V$-$\mathrm{A}$-module which is of finite type as $\mathrm{H}^{*}V$-module such that $dth_{V}M = k+1$. By induction hypothesis, $M$ is isomorphic, as a $DV$-module, to: $\mathbb{F}_{2}[c_1,...,c_k] \otimes \overline{M}^{(c_1,...,c_k)}$. Now, according to \cite{BouZ}, since $dth_{H^{*}V}M = k+1$, the sequence $\{c_1,...,c_{k+1}\}$ is regular for $M$ that means the multiplication by $c_{k+1}$ is injective on the $\mathbb{F}_{2}[c_{k+1},...,c_n]$-module $\overline{M}^{(c_1,...,c_k)}$. The method used for the case "$k=1$" shows that we have an isomorphism of $\mathbb{F}_{2}[c_{k+1},...,c_n]$-modules: $\overline{M}^{(c_1,...,c_k)} \cong \mathbb{F}_{2}[c_{k+1}] \otimes \overline{\overline{M}^{(c_1,...,c_k)}}^{(c_{k+1})}$. \\
So we deduce the isomorphism of $DV$-modules: $M \cong \mathbb{F}_{2}[c_1,...,c_{k+1}] \otimes \overline{M}^{(c_1,...,c_{k+1})}$.
\end{proof}
\begin{re3} The structure of unstable $A$-module in the lemma $2.3.2.3$ is not necessary as it was pointed for us by Professor Jean Lannes in a private communication. One can verify that if $M$ is an $\mathrm{H}^{*}V$-module and if $\{x_{1},...,x_{k}\}$ is a regular sequence for $M$ then $M$ is a free as an $\mathbb{F}_{2}[x_{1},...,x_{k}]$-module.
\end{re3}
\subsection{Technical result}
 Let $V$ be an elementary abelian $2$-group of rank $n \geq 1$, $X$ a finite $V$-CW-complex and $W$ a subgroup of $V$ of codimension one. In this paragraph, we further specify the following notations:
\medskip
\begin{itemize}
  \item $H^{*}V \cong \mathbb{F}_{2}[t_{1},...,t_{n}]$
\medskip
  \item $H^{*}(V/W) \cong \mathbb{F}_{2}[t_{n}]$, $\overline{H^{*}_{V}X}^{H^{*}(V/W)} = \overline{H^{*}_{V}X}^{t_{n}}$ and $\mathcal{T}^{H^{*}(V/W)}\big(H^{*}_{V}X \big)= \mathcal{T}^{t_{n}}\big(H^{*}_{V}X \big)$.
\medskip
  \item $DV \cong \mathbb{F}_{2}[c_{V,1},...,c_{V,n}]$ and $DW \cong \mathbb{F}_{2}[c_{W,1},...,c_{W,n-1}]$
\end{itemize}
\begin{pr2} Let $V$ be an elementary abelian $2$-group of rank $n \geq 1$, let $k \in \mathbb{N}, 1 \leq k \leq n-1$, let $X$ be a finite $V$-CW-complex and $W$ a subgroup of $V$ of codimension one. The following properties are equivalent:
\medskip
\begin{itemize}
  \item [(i)] $dth_{H^{*}W}H^{*}_{W}X \geq k$
\medskip
  \item [(ii)] $\left\{
               \begin{array}{ll}
                 (ii)-a: \; dth_{H^{*}W}\overline{H^{*}_{V}X}^{H^{*}(V/W)} \geq k, & \hbox{}
\medskip\\
                 (ii)-b: \; dth_{H^{*}W}\mathcal{T}^{H^{*}(V/W)}\big(H^{*}_{V}X \big) \geq k. & \hbox{}
               \end{array}
             \right.$
\end{itemize}
\end{pr2}
\medskip
By lemma 2.3.2.3, the proposition 2.4.1 is equivalent to:
\begin{pr2} Let $V$ be an elementary abelian $2$-group of rank $n \geq 1$, let $k \in \mathbb{N}, 1 \leq k \leq n-1$, let $X$ be a finite $V$-CW-complex and $W$ a subgroup of $V$ of codimension one. The following properties, where isomorphisms are isomorphisms of $DW$-modules, are equivalent:
\medskip
\begin{itemize}
\item [(i)] $H^*_{W}X \cong \mathbb{F}_{2}[c_{W,1},...,c_{W,k}] \otimes \overline{H^*_{W}X}^{(c_{W,1},...,c_{W,k})}$.
\medskip
\item [(ii)] $\left\{
               \begin{array}{ll}
                 (ii)-a: \; \overline{H^*_{V}X}^{H^{*}(V/W)} \cong \mathbb{F}_{2}[c_{W,1},...,c_{W,k}] \otimes \overline{\overline{H^*_{W}X}^{H^{*}(V/W)}}^{(c_{W,1},...,c_{W,k})}, & \hbox{}
\medskip\\
                 (ii)-b: \mathcal{T}^{H^{*}(V/W)}(H^*_{V}X) \cong \mathbb{F}_{2}[c_{W,1},...,c_{W,k}] \otimes \overline{\mathcal{T}^{H^{*}(V/W)}(H^*_{V}X)}^{(c_{W,1},...,c_{W,k})}. & \hbox{}
               \end{array}
             \right.$
\end{itemize}
\end{pr2}
\bigskip
Let $V$ be an elementary abelian $2$-group, $X$ a finite $V$-CW-complex and $W$ a subgroup of $V$ of codimension one. Before proving proposition 2.4.2, we introduce the exact Gysin sequence associated to the subgroup $W$ of $V$ of codimension one.
\medskip\\
The inclusion $i: W \hookrightarrow V$ induces a map, denoted $i: X_{hW} \rightarrow X_{hV}$ which is a two-sheet covering:
 $$ V/W \cong \mathbb{Z}/2\mathbb{Z} \rightarrow X_{hW} \overset{i} \rightarrow X_{hV}$$
whose classifying map is $B(\pi) \circ \pi_{V}:  X_{hV} \rightarrow BV \rightarrow B(V/W)$ where $\pi: V \rightarrow V/W$ is the natural projection.\\
Since $H^{*}(V/W) \cong \mathbb{F}_{2}[t_{n}]$, we still denote by $t_{n}$  the non trivial element $(B(\pi) \circ \pi_{V})^{*}(t_{n}) \in H^{1}_{V}X$. The Gysin sequence associated to the precedent covering, is the following sequence of $\mathrm{H}^{*}V$-modules:
\begin{center}
$G(W, V):\; \xymatrix{...\ar[r] & H^{*-1}_{V}X \ar[r]^{t_{n}.} & H^{*}_{V}X \ar[r]^{i^{*}} & H^{*}_{W}X \ar[r]^{tr} & H^{*}_{V}X \ar[r]^{t_{n}.} &...}$
\end{center}
where $tr$ is the transfert (\cite{Sp}, \cite{Z}) and $t_{n}.x = (B(\pi) \circ \pi_{V})^{*}(t_{n}) \smile x$, $x \in H^{*}_{V}X$.
\medskip\\
The exact Gysin sequence $G(W,V)$ induces then a short exact sequence of $\mathrm{H}^{*}W$-modules:
$$ \overline{G}(W,V):\;\; \xymatrix{0 \ar[r] &\overline{H^{*}_{V}X}^{H^{*}(V/W) }
\ar[rr]^-{i^{*}}&&H^{*}_{W}X\ar[rr]^-{tr} && \mathcal{T}^{H^{*}(V/W) }(H^*_{V}X)\ar[r] & 0}
$$
{\bf Proof of proposition 2.4.2}. Let $V$ be an elementary abelian 2-group of rank $n$, $X$ a finite $V$-CW-complex and $W$ a subgroup of $V$ of codimension one. Let $ \overline{G}(W,V)$ the short exact Gysin sequence considered as an exact sequence of $DW$-modules.
$$ \overline{G}(W,V):\; \xymatrix{0 \ar[r] &\overline{H^{*}_{V}X}^{t_{n}} \ar[r]^-{i^{*}}&H^{*}_{W}X \ar[r]^-{tr}  & \mathcal{T}^{t_{n}}(H^*_{V}X) \ar[r] & 0 }$$
{\bf The implication $(ii) \Rightarrow (i)$}. From the depth property 2.2.3 we have: $$dth_{H^{*}W}H^{*}_{W}X \geq min \big(dth_{H^{*}W}\overline{H^{*}_{V}X}^{t_{n}},\; dth_{H^{*}W}\mathcal{T}^{t_{n}}(H^{*}_{V}X)\big)$$
By assumption $(ii)$,  this minimum is $ \geq k$ so $dth_{H^{*}W}H^{*}_{W}X \geq k$. The lemma 2.3.2.3 shows that we have an isomorphism of $DW$-modules: $H^{*}_{W}X \cong \mathbb{F}_{2}[c_{W,1},...,c_{W,k}] \otimes \overline{H^*_{W}X}^{(c_{W,1},...,c_{W,k})}$.
\medskip\\
{\bf The implication $(i) \Rightarrow (ii)$}. By hypothesis we have this isomorphism of $DW$-modules: $H^*_{W}X \cong \mathbb{F}_{2}[c_{W,1},...,c_{W,k}] \otimes \overline{H^*_{W}X}^{(c_{W,1},...,c_{W,k})}$. Considering the sequence "$\mathbb{F}_{2} \otimes _{\mathbb{F}_{2}[c_{W,1}]} \overline{G}(W, V)$", we get:
\begin{itemize}
\medskip
 \item [2.4.2.1] \;\; the following exact seqence of $\mathbb{F}_{2}[c_{W,2},...,c_{W,n-1}]$-modules:
 \clearpage
 $\xymatrix{0 \ar[r] & Tor_{1}^{\mathbb{F}_{2}[c_{W,1}]}(\mathbb{F}_{2}, \mathcal{T}^{t_{n}}(H^*_{V}X)) \ar[r] & \overline{\overline{H^{*}_{V}X}^{t_{n}}}^{(c_{W,1})} \ar[r] &  }$ \\
$\xymatrix{ \mathbb{F}_{2}[c_{W,2},...,c_{W,k}] \otimes \overline{H^*_{W}X}^{(c_{W,1},...,c_{W,k})} \ar[r]  &  \overline{\mathcal{T}^{t_{n}}(H^*_{V}X)}^{(c_{W,1})} \ar[r] & 0 }$
\medskip\\
   and
\medskip
\item [2.4.2.2] \;\;  $Tor_{1}^{\mathbb{F}_{2}[c_{W,1}]}\big(\mathbb{F}_{2}, \overline{H^{*}_{V}X}^{t_{n}} \big)   =  0$
\end{itemize}
\medskip
Since $Tor_{1}^{\mathbb{F}_{2}[c_{W,1}]}\big(\mathbb{F}_{2}, \overline{H^{*}_{V}X}^{t_{n}} \big) \cong \Sigma^{2^{n-2}} \mathcal{T}^{(c_{W,1})}\big(\overline{H^{*}_{V}X}^{t_{n}}\big)$, the point $2.4.2.2$ shows that the multiplication by $c_{W,1}$ on $\overline{H^{*}_{V}X}^{t_{n}}$ is injective which means that: $dth_{H^{*}W}\overline{H^{*}_{V}X}^{t_{n}} \geq 1$. The lemma 2.3.2.3 shows that we have an isomorphism of $DW$-modules:
$$(I1): \; \overline{H^{*}_{V}X}^{t_{n}} \cong \mathbb{F}_{2}[c_{W,1}] \otimes \overline{\overline{H^{*}_{V}X}^{t_{n}}}^{(c_{W,1})}$$
\medskip
Considering "$ \mathbb{F}_{2}[c_{W,1}] \otimes (2.4.2.1)$" and isomorphisme $(I1)$, we get:
$$\begin{array}{lll}
                  Tor_{1}^{\mathbb{F}_{2}[c_{W,1}]}\big(\mathbb{F}_{2}, \mathcal{T}^{t_{n}}(H^*_{V}X) \big)  & = & 0
\medskip\\
 &=& \Sigma^{2^{n-2}} \mathcal{T}^{c_{W,1}}\big( \mathcal{T}^{t_{n}}(H^*_{V}X) \big).
 \end{array}$$
Since $\mathcal{T}^{c_{W,1}}\big( \mathcal{T}^{t_{n}}(H^*_{V}X) \big) = 0$ then $dth_{H^{*}W}\mathcal{T}^{t_{n}}(H^*_{V}X) \geq 1$.  The lemma 2.3.2.3 shows an isomorphism of $DW$-modules:
$$ (I2):\; \mathcal{T}^{t_{n}}\big(H^{*}_{V}X \big) \cong \mathbb{F}_{2}[c_{W,1}] \otimes \overline{\mathcal{T}^{t_{n}}\big(H^{*}_{V}X \big)}^{(c_{W,1})}$$
Considering "$ \mathbb{F}_{2} \otimes_{\mathbb{F}_{2}[c_{W,1}]} \overline{G}(W, V)$"  and the previous isomorphisms $(I1)$ et $(I2)$, we obtain the short exact sequence of $\mathbb{F}_{2}[c_{W,2},...,c_{W,n-1}]$-modules:
$$\xymatrix{0 \ar[r] & \overline{\overline{H^{*}_{V}X}^{t_{n}}}^{(c_{W,1})} \ar[r] & \mathbb{F}_{2}[c_{W,2},...,c_{W,k}] \otimes \overline{H^*_{W}X}^{(c_{W,1},...,c_{W,k})} \ar[r]  &  \overline{\mathcal{T}^{t_{n}}\big(H^{*}_{V}X \big)}^{(c_{W,1})}\ar[r] & 0 }$$
The same method applied to this sequence and from step to step, one obtain the requested isomorphisms of $DW$-modules.
\section{The main result}
Let $V$ be an elementary abelian $2$-group and $X$ a finite $V$-CW-complex. The modulo~2 equivariant cohomology of $X$, is naturally an unstable $H^*V$-$A$-module which is of finite type as $\mathrm{H}^{*}V$-module; let's denote
$dth_{H^*V}H^*_{V}X$ the depth of $H^*_{V}X$ relatively to the augmentation ideal, $\widetilde{H^{*}}V$, de $\mathrm{H}^{*}V$.\\
Let $W$ be a subgroup of $V$, the finite $V$-CW-complex $X$ is also a finite $W$-CW-complexe; thus we have two integers namely $dth_{H^*V}H^*_{V}X$ and $dth_{H^{*}W}H^*_{W}X$. Our main result is:
\begin{th1}
Let $V$ be an elementary abelian $2$-group and $X$ a finite $V$-CW-complex.  Then, for any subgroup $W$ of $V$, we have the inequality:
 $dth_{H^{*}W} H^{*}_{W}X \leq dth_{H^{*}V} H^{*}_{V}X $.
\end{th1}
\subsection{Proof of theorem 3.1}
Let $V$ be an elementary abelian 2-group of rank $n \geq 1$, $X$ a finite $V$-CW-complex and $W$ a subgroup of $V$. We want to show the inequality: $dth_{H^{*}W} H^{*}_{W}X  \leq dth_{H^{*}V} H^{*}_{V}X  $.
\medskip\\
Note that if $dth_{H^{*}W}H^{*}_{W}X = 0$, then  the inequality is trivial. Therefore, we assume \\$dth_{H^*W}H^{*}_{W}X  \geq 1$.
\medskip\\
It is clear that we easily return to the case $W$ is a subgroup of $V$ of codimension one.
\medskip\\
We distinguish two cases, depending on whether, the trivial part $\mathcal{T}^{H^{*}(V/W)} \big(H^{*}_{V}X \big) $ is zero or not.
\subsubsection{Case where $\mathcal{T}^{H^{*}(V/W)} \big(H^{*}_{V}X \big) = 0$}
\paragraph{ } Consider the exact Gysin sequence  $\overline{G}(W, V)$ of $H^{*}W$-modules:
$$ \overline{G}(W, V):\;\; \xymatrix{0 \ar[r] & \overline{H^{*}_{V}X }^{H^{*}(V/W)} \ar[rr]^-{i^{*}}&&H^{*}_{W}X \ar[rr]^-{tr} && \mathcal{T}^{H^{*}(V/W)} \big( H^{*}_{V}X  \big) \ar[r] & 0}$$
Since $\mathcal{T}^{H^{*}(V/W)} \big(H^{*}_{V}X \big) = 0$, this exact Gysin sequence yields an isomorphism of $H^{*}W$-modules: $\overline{H^{*}_{V}X }^{H^{*}(V/W)} \cong H^{*}_{W}X $. So we have:
$$\begin{array}{lll}
                    dth_{H^{*}W}H^{*}_{W}X  & = & dth_{H^{*}W} \big( \overline{H^{*}_{V}X }^{H^{*}(V/W)} \big)
\medskip\\
                   & = & dth_{H^{*}V}H^{*}_{V}X  - 1 \;\; \text{(see proposition 2.3.1.1).}
\end{array}$$
\medskip
This shows that $dth_{H^{*}V}H^{*}_{V}X = dth_{H^{*}W}H^{*}_{W}X  + 1 > dth_{H^{*}W}H^{*}_{W}X $.
\subsubsection{Case where $\mathcal{T}^{H^{*}(V/W)} \big(H^{*}_{V}X \big) \neq 0$}
 \paragraph{ }  Let $V$ be an elementary abelian 2-group of rank $n \geq 1$, $X$ a finite $V$-CW-complex and $W$ a subgroup of $V$ of codimension one. Recall the notations:
\medskip
\begin{itemize}
  \item $H^{*}V \cong \mathbb{F}_{2}[t_{1},...,t_{n}]$.
\medskip
  \item $H^{*}(V/W) \cong \mathbb{F}_{2}[t_{n}]$.
\medskip
  \item Dickson algebras: $DV \cong \mathbb{F}_{2}[c_{V,1},...,c_{V,n}]$ and $DW \cong \mathbb{F}_{2}[c_{W,1},...,c_{W,n-1}]$.
\end{itemize}
\medskip
Before continuing the proof of theorem 3.1, we need to introduce the following graded polynomial algebra: $ \widetilde{D}V = \mathbb{F}_{2}[c_{W,1},...,c_{W,n-1}, t_{n}] $
\medskip\\
We have these properties:
\medskip\\
(3.1.2-a) Let $s: V/W \hookrightarrow V$ be a linear section  of the natural projection $q: V \twoheadrightarrow V/W$ and $GL(V)$ the groupe of automorphismes of $V$. let's write $V \cong W \oplus s(V/W)$ and denote $GL\big(W, s(V/W)\big)$ the subgroup of $GL(V)$ image of the homomorphism
$$ GL(W) \rightarrow GL(V\cong W \oplus s(V/W)),\; B \mapsto \left(
                                       \begin{array}{cc}
                                         B & 0 \\
                                         0 & 1 \\
                                       \end{array}
                                     \right).$$
We have: $\widetilde{D}V = (H^{*}V)^{GL\big(W, s(V/W)\big)} \cong DW \otimes H^{*}(V/W)$. It results:
\medskip\\
(i) The following inclusions: $DV \hookrightarrow \widetilde{D}V \hookrightarrow H^{*}V $.
\medskip\\
(ii) $H^{*}V$ is a free $\widetilde{D}V$-module of finite type. Since,
$\widetilde{D}V \cong DW \otimes H^{*}(V/W)$, $H^{*}V \cong H^{*}W \otimes H^{*}(V/W)$ and $H^{*}W$ is a free $DW$-module of finite type, (voir \cite{Se1}, \cite{R}).
\medskip\\
(iii) The natural inclusion $i: W \hookrightarrow V$ induces in  modulo $2$ cohomology a morphism of algebras $i^{*}_{\mid}: \widetilde{D}V \rightarrow DW$ and we have the following commutative diagram:
\begin{center}
$\xymatrix{ H^{*}(V/W) \, \ar@{^{(}->}[r]  \, \ar[d]^{=}&  \widetilde{D}V \ar@{^{(}->}[d] \ar@{->>}[rr]^-{i^{*}_{\mid}}&& DW \ar@{^{(}->}[d] \\
H^{*}(V/W) \, \ar@{^{(}->}[r]^-{q^{*}}& H^{*}V \ar@{->>}[rr]^-{i^{*}}&& H^{*}W }$
\end{center}
\medskip
In particular, for any $\widetilde{D}V$-module $N$, we have an isomorphism of $DW$-modules:\\
 $DW \otimes _{\widetilde{D}V} N \cong \mathbb{F}_{2} \otimes _{H^{*}V/W} N$.
\medskip\\
 (3.1.2-b) Let $M$ be an $\mathrm{H}^{*}V$-module  of finite type, then $M$ is a $\widetilde{D}V$-module of finite type and we have, according to \cite{Se 2}, $dth_{H^{*}V}M = dth_{\widetilde{D}V}M$.
\medskip\\
Let's return to the proof of theorem 3.1 in the case where $\mathcal{T}^{H^{*}(V/W)} \big(H^{*}_{V}X \big) \neq 0$. Let $k \in \mathbb{N}, \; 1 \leq k \leq rg(W)$, and assume that $dth_{H^*W}H^{*}_{W}X = k$, we will show that $dth_{H^{*}V}H^{*}_{V}X \geq k$ or equivalently, $dth_{\widetilde{D}V}H^{*}_{V}X \geq k$, where $\widetilde{D}V$ denoting the algebra introduced above.
\medskip\\
According to Auslander-Buchsbaum formula (see 2.2.1): $dth_{\widetilde{D}V}H^{*}_{V}X  +  dp_{\widetilde{D}V}H^{*}_{V}X = n$, we will show that: $dp_{\widetilde{D}V}H^{*}_{V}X = n -  dth_{\widetilde{D}V}H^{*}_{V}X \leq n - k$. Equivalently, we will prove that: $Tor_{i}^{\widetilde{D}V}(\mathbb{F}_{2}, H^{*}_{V}X) = 0 \; \mbox{for} \; i \geq n - k + 1 $.
\medskip\\
Since $dth_{H^*W}H^{*}_{W}X = k$, $1 \leq k \leq n-1$, then by 2.3.2.3 and 2.4.2, we have the following isomorphisms of $DW$-modules:
 \medskip
\begin{itemize}
  \item [3.1.2-c] $H^{*}_{W}X \cong \mathbb{F}_{2}[c_{W,1},...,c_{W,k}] \otimes \overline{ H^{*}_{W}X }^{(c_{W,1},...,c_{W,k})}$.
 \medskip
  \item  [3.1.2-d] $\overline{ H^{*}_{V}X }^{ H^{*}(V/W) } \cong  \mathbb{F}_{2}[c_{W,1},...,c_{W,k}] \otimes
\overline{ \overline{ H^{*}_{V}X }^{ H^{*}(V/W) } }^{ (c_{W,1},...,c_{W,k}) }$.
\medskip
  \item [3.1.2-e] $\mathcal{T}^{H^{*}(V/W)} \big( H^{*}_{V}X \big) \cong \mathbb{F}_{2}[c_{W,1},...,c_{W,k}] \otimes \overline{ \mathcal{T}^{H^{*}(V/W)} \big( H^{*}_{V}X \big) }^{(c_{W,1},...,c_{W,k})}$.
\end{itemize}
\medskip
To show that: $Tor_{i}^{\widetilde{D}V}(\mathbb{F}_{2}, H^{*}_{V}X) = 0 \; \mbox{for} \; i \geq n - k + 1 $, we will consider a spectral sequence (so-called Grothendick spectral sequence) of change of algebras  (see \cite{CE}).
\medskip\\
 Let $V$ be an elementary abelian $2$-group of rank $n \geq 1$, $X$ a finite $V$-CW-complex and $W$ a subgroup of $V$ of codimension one. The natural inclusion $W \hookrightarrow V$ induces in modulo $2$ cohomology an homomorphism of algebras $\widetilde{D}V \rightarrow DW$ (see 3.1.2-a). We have:
\begin{pr3} (\cite{CE}) With the previous notations, there exists a spectral sequence of the first quadrant with term
$E_{p, q}^{2} = Tor_{p}^{DW} \big( \mathbb{F}_{2},  Tor_{q}^{\widetilde{D}V} \big(DW, H^{*}_{V}X \big) \big)$
and converges to $Tor_{p+q}^{\widetilde{D}V}\big(\mathbb{F}_{2}, H^{*}_{V}X \big)$.
\end{pr3}
\medskip
The next result is useful for computing the $E^{2}$-term of the Grothendick spectral sequence.
\begin{le3}
$$Tor_{q}^{\widetilde{D}V} \big(DW, H^{*}_{V}X \big) = \left\{
                                                \begin{array}{ll}
 \overline{H^{*}_{V}X}^{t_{n}}, & \hbox{si} \; q=0,
\medskip\\
\Sigma \mathcal{T}^{t_{n}} \big( H^{*}_{V}X \big), & \hbox{si}\; q=1,
\medskip\\
                                                  0, & \hbox{si}\; q \geq 2.
                                                \end{array}
                                              \right.$$
\end{le3}
\begin{proof} By definition, we have:
\medskip\\
$\begin{array}{lll}

                    Tor_{q}^{\widetilde{D}V} \big(DW, H^{*}_{V}X \big) & = & H_{q}\big(DW \otimes _{\widetilde{D}V} L_{\bullet}\big) \; \text{where} \; L_{\bullet} \twoheadrightarrow H^{*}_{V}X \; \text{is a}\;  \widetilde{D}V  \text{-free resolution}\\
  & &
                    \text{of} \; H^{*}_{V}X.
\medskip\\
& =& H_{q}\big(\mathbb{F}_{2}  \otimes _{H^{*}(V/W)} L_{\bullet}\big) \; \text{because}\; DW \otimes _{\widetilde{D}V} L_{\bullet} \cong \mathbb{F}_{2} \otimes _{H^{*}V/W} L_{\bullet}.
\medskip\\
& \cong &  Tor_{q}^{H^{*}(V/W)} \big(\mathbb{F}_{2}, H^{*}_{V}X \big) \; \text{because}\; \widetilde{D}V \;  \text{is an}\;  H^{*}(V/W)  \text{-free module}
\end{array}$
\medskip\\
Since $H^{*}(V/W) \cong \mathbb{F}_{2}[t_{n}]$, we have:
\medskip\\
$Tor_{q}^{H^{*}(V/W)} \big(\mathbb{F}_{2}, H^{*}_{V}X \big)$=$\left\{
                                                              \begin{array}{ll}
                                                                \overline{H^{*}_{V}X}^{t_{n}}, & \hbox{if} \; q=0.
\medskip\\
                                                               \Sigma \mathcal{T}^{t_{n}} \big( H^{*}_{V}X \big), & \hbox{if} \; q=1.
\medskip\\
                                                                0, & \hbox{if}\; q \geq 2.
                                                              \end{array}
                                                            \right.$
\end{proof}
\medskip
It follows that $E_{2}^{p,q} = 0$ for $q \ge 2$. Thus, we have the exact sequence (see \cite[p. 328]{CE}):
\smallskip
$$ (S):\;\; \xymatrix{ \cdots \ar[r] & E_{2}^{*-2,1} \ar[r] & E_{2}^{*,0} \ar[r] & Tor_{*}^{\widetilde{D}V}(\mathbb{F}_{2}, H^{*}_{V}X) \ar[r] & E_{2}^{*-1,1} \ar[r] & \cdots }
$$
\smallskip
The isomorphisms of $DW$-modules $3.1.2$-d and $3.1.2$-e show that we have the following ones:
\bigskip\\
$\begin{array}{lll}
                    E_{2}^{*,0}  & \cong & Tor_{*}^{DW}\big(\mathbb{F}_{2},  \mathbb{F}_{2}[c_{W,1},...,c_{W,k}] \otimes \overline{\overline{ H^{*}_{V}X }^{ t_{n} } }^{(c_{W,1},...,c_{W,k})}\big)
\medskip\\
                   & \cong & Tor_{*}^{\mathbb{F}_{2}[c_{W, k+1},...,c_{W, n-1}]}\big(\mathbb{F}_{2},  \overline{\overline{ H^{*}_{V}X }^{ t_{n} } }^{(c_{W,1},...,c_{W,k})} \big).
\end{array}$
\bigskip\\
$\begin{array}{lll}
E_{2}^{*-1,1}& \cong & Tor_{*-1}^{DW}\big( \mathbb{F}_{2},  \mathbb{F}_{2}[c_{W,1},...,c_{W,k}] \otimes  \overline{\mathcal{T}^{t_{n}}(H^{*}_{V}X)}^{(c_{W,1},...,c_{W,k})} \big)
\medskip\\
                   & \cong & Tor_{*-1}^{\mathbb{F}_{2}[c_{W, k+1},...,c_{W, n-1}]}\big(\mathbb{F}_{2},  \overline{\mathcal{T}^{t_{n}}(H^{*}_{V}X)}^{(c_{W,1},...,c_{W,k})} \big).
\end{array}$
\medskip\\
This proves that: $\left\{
                 \begin{array}{ll}
                   E_{2}^{*,0} = 0, & \hbox{for $* \geq n-k$}
                   \medskip\\
                   E_{2}^{*-1,1} = 0, & \hbox{for $* \geq n-k+1$}
                 \end{array}
               \right.$
\medskip\\
The exact sequence $(S)$ then shows that $ Tor_{*}^{\widetilde{D}V}(\mathbb{F}_{2}, H^{*}_{V}X) = 0$ for $* \geq n-k+1$, which is equivalent to $dp_{\widetilde{D}V}H^{*}_{V}X \leq n-k$ or to $dth_{\widetilde{D}V}H^{*}_{V}X = dth_{H^{*}V}H^{*}_{V}X \geq k$.
\section{Some remarks and consequences}
Let $V$ be an elementary abelian $2$-group, $X$ a finite $V$-CW-complex and $W$ a subgroup of $V$. We remark a certain analogy between the propreties of "topological" action of $V$ on $X$ and the depth of the $H^{*}V$-module of finite type $H^{*}_{V}X$.
\paragraph{\bf Trivial action} The analogue of the trivial action of $V$ on $X$ is an action of $V$ on $X$ such that $H^{*}_{V}X$ is Cohen-Macaulay. In fact, we know that, if $V$ acts trivially on $X$, then any subgroup $W$ of $V$ acts trivially on $X$. We have:
\begin{le1} Let $V$ be an elementary abelian $2$-group and $X$ a finite $V$-CW-complex. If the $H^{*}V$-module of finte type $H^{*}_{V}X$ is Cohen-Macaulay, then it is the same for the $H^{*}W$-module of finite type $H^{*}_{W}X$, for any subgroup $W$ of $V$.
\end{le1}
\begin{proof} We have just to prove the lemma for a subgroup $W$ of $V$ of codimension one. Since the $H^{*}V$-module of finite type $H^{*}_{V}X$ is Cohen-Macaulay, we have: $dth_{H^{*}V}H^{*}_{V} = rk(V)$ or equivalently $dp_{H^{*}V}H^{*}_{V}X = 0$, by Auslander-Buchsbaum's formula (see 2.2.1)). Hence $H^{*}_{V}X$ is a free $H^{*}V$-module so $\mathcal{T}^{H^{*}(V/W)} \big(H^{*}_{V}X \big) = 0 $. The short Gysin exact sequence, $\overline{G}(W, V)$, gives the isomorphism of $H^{*}W$-modules: $\overline{H^{*}_{V}X }^{H^{*}(V/W)} \cong H^{*}_{W}X $. Therefore, we deduce:
$$\begin{array}{lll}
                    dth_{H^{*}W}H^{*}_{W}X  & = & dth_{H^{*}W} \big( \overline{H^{*}_{V}X }^{H^{*}(V/W)} \big)
\medskip\\
                   & = & dth_{H^{*}V}H^{*}_{V}X  - 1 \;\; \text{(voir proposition 2.3.1.1).}
\medskip\\
                   & = & rk(V) - 1.
\end{array}$$
\medskip
which proves that $dth_{H^{*}W}H^{*}_{W}X = rk(V)-1=rk(W)$ that is the $H^{*}W$-module  $H^{*}_{W}X$ is Cohen-Macaulay.
\end{proof}
\medskip
\paragraph{\bf Free action} The analogue of free action of $V$ on $X$ is an action of $V$ on $X$ such that $dth_{H^{*}V}H^{*}_{V}X = 0$. In fact, we know that, if $V$ acts freely on $X$, then any subgroup $W$ of $V$ acts freely on $X$. By theorem 3.1, we have the following result:
\begin{le1} Let $V$ be an elementary abelian $2$-group and $X$ a finite $V$-CW-complex. If $dth_{H^{*}V}H^{*}_{V}X = 0$, then $dth_{H^{*}W}H^{*}_{W}X = 0$, for any subgroup $W$ of $V$.
\end{le1}
\medskip
More generally, let's define the following:
\begin{de1} Let $V$ be an elementary abelian $2$-group, $X$ a finite $V$-CW-complex and $k \in \mathbb{N}, \; 0 \leq k \leq rk(V)$. We say that the action of $V$ on $X$ is a {\bf $k$-action} if $dth_{H^{*}V}H^{*}_{V}X \leq k$.
\end{de1}
\medskip
A consequence of the theorem 3.1 is:
\begin{pr1} Let $V$ be an elementary abelian $2$-group, $X$ a finite $V$-CW-complex and $k \in \mathbb{N}, \; 0 \leq k \leq rk(V)$.  If the action of $V$ on $X$ is a $k$-action, then the same is true for any subgroup of $V$.
\end{pr1}

\end{document}